\documentclass[8pt]{article}
\usepackage{indentfirst,latexsym,bm}
\usepackage{amsfonts}
\usepackage{amssymb}
\usepackage[leqno]{amsmath}
\usepackage{dsfont}
\usepackage{amsthm}
\usepackage{times}
\usepackage[all]{xy}
\usepackage{hyperref}
\usepackage{color,xcolor}
\newtheorem{theorem}{Theorem}[section]
\newtheorem{corollary}[theorem]{Corollary}
\newtheorem{lemma}[theorem]{Lemma}
\newtheorem{proposition}[theorem]{Proposition}
\theoremstyle{definition}
\newtheorem{remark}[theorem]{Remark}
\newcommand{\C}{\mathbb C}
\newcommand{\cC}{\mathcal C}
\newcommand{\cD}{\mathcal D}
\newcommand{\cE}{\mathcal E}
\newcommand{\cO}{\mathcal O}
\newcommand{\sVec}{\mathrm{sVec}}
\newcommand{\co}{\operatorname{co}}
\newcommand{\id}{\mathrm{id}}
\newcommand{\eps}{\varepsilon}
\newcommand{\cop}{\mathrm{cop}}

\DeclareMathOperator{\FPdim}{FPdim}
\DeclareMathOperator{\Rep}{Rep}
\DeclareMathOperator{\Vect}{Vec}
\DeclareMathOperator{\Tr}{Tr}
\DeclareMathOperator{\im}{im}
\DeclareMathOperator{\ord}{ord}
\DeclareMathOperator{\lcm}{lcm}

\numberwithin{equation}{section}
\begin{document}
\title{On quasi-triangular Hopf algebras of square-free dimensions}
\date{}
\author{Jiacheng Bao, Zhiqiang Yu, and  Ying Zheng}
\maketitle
\abstract
Let $H$ be a quasi-triangular Hopf algebra of sqaure-free dimension over $\C$. In this note, we show that $H$ is semisimple and  cocommutative. Equivalently, a braided finite tensor category $\cC$ of square-free Frobenius--Perron dimension with a fiber functor is
tensor equivalent to the representation category of a finite group.


\section{Introduction}
Let $H$ be a finite-dimensional Hopf algebra over $\C$. Then the category $\Rep(H)$ of finite-dimensional $H$-modules is a tensor category in sense of \cite{EGNO}. Moreover, one can recover the structure of $H$ by the Tannakian  reconstruction. Indeed, there is a bijective correspondence between finite-dimensional Hopf algebras and finite tensor categories with a fiber functor, i.e. a faithful exact tensor functor to the category $\Vect$ of finite-dimensional vector spaces, see \cite{EGNO}.

It was asked by N. Andruskiewitsch that whether a Hopf algebra $H$ of square-free dimension is always semisimple \cite[Question 6.5]{Andrusk}. Up to now,  only partial affirmative answer has   been obtained in \cite{Ng2p,Ngpq,Zhu}. If we further assume that $H$  is an odd-dimensional quasi-triangular Hopf algebra, then $H$ is semisimple and cocommutative \cite{Natale} (see also \cite{WZ}), that is, $H\cong\C[G]$ for a finite group $G$. Recall that a quasi-triangular Hopf algebra is a Hopf algebra with an $R$-matrix, which was introduced by Drinfeld \cite{DrinfeldICM}, see \cite{Radford} and the references therein for further details on these Hopf algebras. Notice that the representation category of finite-dimensional quasi-triangular Hopf algebras are braided tensor categories \cite{DGNO,EGNO}.
Conversely, a braided finite tensor
category together with a fiber functor to $\Vect$ reconstructs a finite-dimensional
quasi-triangular Hopf algebra \cite{EGNO}. Here a fiber functor means a faithful exact tensor functor.

In this note, combined conclusions of \cite{Ng2p,WZ} with the theory of braided tensor   categories \cite{DGNO,EGNO},  we  remove the oddness assumption from
 \cite[Theorem 1.2]{Natale}.  Namely,
\begin{theorem}
Let $(H,R)$ be a quasitriangular Hopf algebra over $\C$. If $\dim (H)$ is
square-free, then $H$ is semisimple and cocommutative. Therefore,
$
H\cong\C [G]$
as a Hopf algebra, for a finite group $G$ . Consequently, every braided finite tensor category $\cC$ over $\C$
of square-free Frobenius--Perron dimension admitting a fiber functor
is a fusion category and is tensor equivalent to $\Rep(G)$ for a
finite group $G$.
\end{theorem}

The proof has two parts. For semisimplicity, we combine
\cite[Proposition~6.3]{WZ} with the antipode parity arguments in
\cite{Ngpq,Ng2p}. These arguments exclude the remaining possibilities
for a minimal counterexample in even square-free dimension. For
cocommutativity, we use de-equivariantization of  a maximal Tannakian subcategory
of the M\"uger center, the resulting   braided fusion category
is non-degenerate or slightly degenerate, and square-freeness forces
it to be pointed. We then obtain a cocentral Hopf algebra extension
with cyclic function-algebra kernel, by  applying
\cite[Lemma~5.8]{Natale11}, we get the cocommutativity.

This note is organized as follows. In Section \ref{section2}, we recall some basic notations of quasi-triangular Hopf algebras, (braided) tensor categories. In Section \ref{section3}, we prove the main result in Theorem \ref{thm:hopf} and Theorem \ref{thm:cocommutative}.

\section{Preliminaries}\label{section2}
In this section, we recall some basic notions of quasi-triangular Hopf algebras and braided finite tensor category over $\C$, the readers are referred to \cite{EGNO,Radford} for details. And all Hopf algebras are assumed to be finite-dimensional.

For a Hopf algebra $H$ over $\C$, we use the notation $\Delta(h)=h_{(1)}\otimes h_{(2)}$, $h\in H$,  $S_H$ and $\varepsilon$ denote  the antipode and counit of $H$, respectively. Moreover, for any Hopf algebra $H$, we use   $H^\cop$ to denote the  cooposite Hopf algebra of $H$, $\Delta^\cop(h)=h_{(1)}\otimes h_{(2)}$. And we use $G(H)$ for its group
of group-like elements, and $Z(H)$ for its algebra center.

Let $(H,R)$ be a quasi-triangular Hopf algebra with  $R$-matrix. Namely, $R=\sum_iR_i\otimes R^i\in H\otimes H$ is an invertible element and satisfies the following conditions:
\begin{equation}\label{eq:qt}
\begin{split}
&(\Delta\otimes\id)(R)=R_{13}R_{23},\quad
(\id\otimes\Delta)(R)=R_{13}R_{12},\\
&(\id\otimes\varepsilon)(R)=(\varepsilon\otimes \id)(R)=1,\quad \Delta^{\cop}(h)=R\Delta(h)R^{-1},\forall h\in H.
\end{split}
\end{equation}
If $R_{21}R=1\otimes 1$,then $H$ is a triangular Hopf algebra, where $R_{21}=\sum R^i\otimes R_i$.
For a Hopf quotient
$\pi:H\twoheadrightarrow B$, namely, $\pi:H\to B$ is a surjective homomorphism of Hopf algebras, we write
$
H^{\co\pi}=\{h\in H|(\id\otimes\pi)\Delta(h)=h\otimes1\}$.
Then $H^{\co \pi}$ is a left coideal subalgebra of $H$.

Let $\cC$ be a finite tensor category, that is, a $\C$-linear finite abelian rigid monoidal category with simple tensor unit. Let $\text{Gr}(\cC)$ be the Grothendieck ring of $\cC$. Then there is a unique homomorphism $\FPdim(-):\text{Gr}(\cC)\to\C$ such that $\FPdim(X)$ is a positive algebraic integer \cite{EGNO}, $X\in\cO(\cC)$, $\FPdim(X)$ is called the Frobenius-Perron dimension of $X$. The Frobenius-Perron dimension of $\cC$
is defined as \begin{align*}
\FPdim(\cC):=\sum_{X\in\cO(\cC)}\FPdim(X)\FPdim(P(X)),
\end{align*}
where $P(X)$ is the projective cover of $X$.  If $\cC=\Rep(H)$, then $\FPdim(X)=\dim_\C(X)$ and $\FPdim(\cC)=\dim(H)$.

A fusion category   $\cC$ is a finite semisimple tensor category. Recall that a finite-dimension Hopf algebra $H$ is semisimple (i.e. isomorphic to a direct sum of simple algebras) if and only if it dual Hopf algebra $H^*$ is semisimple \cite{Radford}, if and only if the category of finite-dimensional left $H$-modules $\Rep(H)$ is a fusion category. Fusion category  $\cC$ is said to be pointed  if all of its simple objects are invertible, equivalently, $\FPdim(X)=1$, $\forall X\in\cO(\cC)$. 

A braided tensor category $\cC$ is a tensor category with  a braiding $c$. For example, let $(H,R)$ be a quasi-triangular Hopf algebra, then  $\Rep(H)$ is a braided tensor category with braiding $c=\tau\circ R$, where $\tau$ is the reflection of vector spaces, namely, \begin{align*}
c_{V,W}(v\otimes w)=\sum_iR^iw\otimes R_iv, \quad v\in V,w\in W.
\end{align*}
For a braided tensor  category $\cC$, its M\"uger center $\cC'$ consists
of objects $X$ for which
\[
c_{Y,X}c_{X,Y}=\id_{X\otimes Y},\qquad\forall Y\in\cC.
\]
The category is non-degenerate when $\cC'=\Vect$ and slightly degenerate
when $\cC'=\sVec$, where $\Vect,\sVec$ are the categories of finite-dimensional vector spaces and super vector spaces over $\C$, respectively. If $\cC=\cC'$, then $\cC$ is symmetric.

\section{Main results}\label{section3}


Let $H$ be a finite-dimensional Hopf algebra. For $a\in H$, let $r(a)$ be right multiplication by $a$ on $H$.
We use the standard trace criterion: if $H$ is not semisimple, then
\begin{equation}\label{eq:trace}
\Tr(S_H^2\circ r(a))=0,\qquad\forall a\in H.
\end{equation}
In particular, $\Tr(S_H^2)=0$, see \cite{LarRad,Radford} for details, also \cite[Theorem~1.2]{Ng2p}.

\begin{lemma}\label{lem:involution}
Suppose that Hopf algebra $H$ is not semisimple and contains a group-like element $g$
of order $2$. Set $e=(1+g)/2$. Then
\begin{equation}\label{eq:idempotent}
\dim (He)=\frac{\dim(H)}{2},\quad
S_H^2(He)=He,\quad
\Tr(S_H^2|_{He})=0.
\end{equation}
If, in addition, $S_H^8=\id_H$, then $4\mid\dim (H)$.
\end{lemma}

\begin{proof}
Let $G:=\langle g\rangle$. Then the Nichols--Zoeller freeness theorem \cite{NZ} says that $H$ is free as a right module over
the Hopf subalgebra $\C [G]$, the group algebra of $G$. Then $2\dim (He)=\dim(H)$. Since $S_H^2(g)=g$, $He$ is stable under $S_H^2$, then $S_H^2(He)=He$ as $S_H$ is bijective.
The operator $r(e)$ is the projection onto $He$ and commutes with
$S_H^2$. Thus \eqref{eq:trace} gives the trace assertion.
Now let $T=S_H^2|_{He}$ and assume $S_H^8=\id_H$. The operator $T$ is
diagonalizable with eigenvalues in $\{1,-1,i,-i\}$. If their multiplicities
are $a,b,c,d$, respectively, then
\[
0=\Tr(T)=(a-b)+i(c-d).
\]
Hence $a=b$ and $c=d$, so $\dim (He)=2(a+c)$ is even. Therefore
$4\mid\dim (H)$.
\end{proof}

\begin{lemma}\label{lem:s4}
If Hopf algebra $H$ is not semisimple and $S_H^4=\id_H$, then $4\mid\dim (H)$.
Equivalently, if $4\nmid\dim(H)$ and $S_H^4=\id_H$, then $H$ is
semisimple.
\end{lemma}

\begin{proof}
Notice that the eigenvalues of $S_H^2$ are $1$ and $-1$, and its trace is zero by assumption.
Consequently, the two eigenspaces have the same dimension.

Let $g\in G(H)$ be the distinguished group-like element.
If $\ord(g)$ is even, then $\langle g\rangle$ contains an involution.
Since $S_H^8=\id_H$, $4\mid \dim(H)$ by Lemma~\ref{lem:involution}. If $g=1$, then $H^*$ is unimodular. So the squared antipode $S_{H^*}^2$ has order $2$, since $H^*$ is not semisimple, $S_H^2$ cannot   be the identity.
\cite[Proposition~2.1]{Ng2p} therefore gives $4\mid\dim (H)^*=\dim(H)$. If $n=\ord(g)>1$ is odd. Recall that the index of
$H$ is the least positive integer $t$ for which both $S_H^{4t}=\id_A$
and $g^t=1$. Here the index for $H$ is $n$. By
\cite[Corollary~2.2]{Ngpq}, the subspace
$V:=\{a\in H|S_H^{2n}(a)=-a\}$
has even dimension. Since $n$ is odd and $S_H^4=\id_H$, we have
$S_H^{2n}=S_H^2$. Thus  $V$ is exactly the eigenspace of $S_H^2$ with eigenvalue $1$,   and $\dim(H)=2\dim(V)$ is divisible by $4$.
\end{proof}


The following proposition is a direct result of \cite[Proposition~6.3]{WZ}.
\begin{proposition}\label{prop:wz}
Let $(H,R)$ be a quasitriangular Hopf algebra such that
\[
G(H)\cap Z(H)=\{1\},\qquad G(H^*)\cap Z(H^*)=\{\eps\}.
\]
Then for every odd prime $p$ dividing $|G(H^*)|$, we have
$p^2\mid\dim (H)$.
\end{proposition}

Now we are ready  to give the semisimplicity of quasi-triangular Hopf algebras of square-free dimensions.
\begin{theorem}\label{thm:hopf}
Every finite-dimensional quasi-triangular Hopf algebra over $\C$ of
square-free dimension $d$ is semisimple.
\end{theorem}

\begin{proof}
If $\dim(H)$ is odd, then this was proved in \cite[Theorem 1.2]{Natale}, see also \cite[Theorem 6.8]{WZ}. Assume $\dim(H)$ is even below.
Suppose that $(H,R)$ is a non-semisimple quasi-triangular Hopf algebra with smallest dimension $d$.

Following \cite{Natale, Radford}, let
$
f_R\colon H^{*\cop}\to H$, $f_R(\varphi)=(\varphi\otimes\id)(R),~\forall \varphi \in H^{*\cop}$.
Let $H_+=\im(f_R)$ and let $H_-$ be the Hopf subalgebra spanned by the
left tensorands of $R$. Then the  minimal quasi-triangular Hopf subalgebra of $H$ is
$
H_R=H_-H_+=H_+H_-$ \cite[Proposition 12.2.13]{Radford},  which is a quotient of the Drinfeld double $D(H_+)$, see \cite{Radford93}. We show that $H_R=H$. Indeed, if $H_R\subsetneq H$, then $\dim (H_R)$ is a proper square-free divisor
of $d$, then the minimality of $\dim(H)$ makes $H_R$ semisimple. Meanwhile, the semisimplicity of $H_R$ implies $S_H^4=\id_H$
by \cite[Theorem 4.4]{Garcia}. Then Lemma~\ref{lem:s4} would then imply
$4\mid d$, it is a contradiction. Therefore $H_R=H$.
Let $r=\dim (H_+)$. Notice that the Nichols--Zoeller freeness theorem \cite{NZ} applied to the subalgebra
$H_+\subseteq H$ and dually to the quotient $D(H_+)\twoheadrightarrow H$,
gives $
r\mid d$ and $d\mid r^2$, hence $r=d$. Thus $f_R:H^{*\cop}
\to H$ is a Hopf algebra isomorphism. In particular, $H^*$ admits a quasi-triangular structure and $
|G(H)|=|G(H^*)|$.

Let $B$ be either $H$ or $H^*$. Both are non-semisimple quasi-triangular
Hopf algebras of dimension $d$. Suppose that $G(B)\cap Z(B)$ contains
a nonidentity element. Taking a suitable power gives a central grouplike
$z$ of prime order $p$. Then $K=\C[G]$ is a normal
semisimple Hopf subalgebra, $G:=\langle z\rangle$, and there is an exact sequence
\[
\C\longrightarrow K\longrightarrow B
\longrightarrow \overline{H}:=B/BK^+\longrightarrow\C,
\qquad \dim( \overline{H})=N/p,
\]
where $K^+=\ker(\eps|_K)$. The quotient $\overline{H}$ inherits a quasi-triangular
structure. Its dimension is a smaller square-free integer, so $\overline{H}$ is
semisimple by minimality of $\dim(H)$. Since an extension of semisimple Hopf algebras
is semisimple, $B$ is
semisimple, a contradiction. Consequently,
\begin{equation}\label{eq:centers}
G(H)\cap Z(H)=\{1\},\qquad
G(H^*)\cap Z(H^*)=\{\eps\}.
\end{equation}

By Proposition~\ref{prop:wz} and \eqref{eq:centers}, any odd prime
$p\mid |G(H^*)|$ would satisfy $p^2\mid d$.
Furthermore, $|G(H^*)|\mid d$ by Nichols--Zoeller freeness theorem \cite{NZ}. Hence
\begin{equation}\label{eq:one-or-two}
|G(H)|=|G(H^*)|\in\{1,2\}.
\end{equation}

Let $g\in G(H)$ and $\alpha\in G(H^*)$ be the distinguished group-like elements.
Then Radford's fourth-power formula \cite[Theorem 10.5.6]{Radford} implies $
\ord(S_H^4)\mid\lcm\bigl(\ord(g),\ord(\alpha)\bigr)$.
If the two groups in \eqref{eq:one-or-two} are trivial, then
$g=1$ and $\alpha=\eps$. Thus $S_H^4=\id_H$ by \cite[Theorem 10.5.6]{Radford}, contradicts to
Lemma~\ref{lem:s4} and square-freeness of $d$. If the two groups have order $2$, then
$S_H^8=\id_H$, however,
Lemma~\ref{lem:involution} yields $4\mid d$,   impossible.
\end{proof}


We recall the following freeness consequence, which is useful because
a coinvariant algebra need not yet be a Hopf subalgebra. If
$\pi:H\twoheadrightarrow B$ is a finite-dimensional Hopf quotient and
$L$ is a Hopf subalgebra contained in $H^{\co\pi}$, then
\begin{equation}\label{eq:coinvariant-div}
\dim (H^{\co\pi})=\frac{\dim (H)}{\dim( B)},
\qquad
\dim (L)\mid\dim( H^{\co\pi}).
\end{equation}
The second statement is the relative Hopf-module form of
Nichols--Zoeller freeness theorem \cite{NZ}, applied to the left $H$-comodule and left
$L$-module $H^{\co\pi}$, see
\cite[Lemma~5.2]{Natale}.

The following lemma is a generalization of \cite[Proposition 4.10]{Natale}.
\begin{lemma}\label{lem:cocentral}
Let $(H,R)$ be a finite-dimensional quasi-triangular Hopf algebra and
$\pi:H\twoheadrightarrow B$ a Hopf quotient. Suppose that
\[
(\pi\otimes\pi)(R)=1\otimes1,\qquad
(\pi\otimes\id)(R_{21}R)=1\otimes1,
\qquad
\gcd\left(\dim (B),\frac{\dim (H)}{\dim( B)}\right)=1.
\]
Then
\begin{equation}\label{eq:mixed-trivial}
(\pi\otimes\id)(R)=1\otimes1=(\id\otimes\pi)(R).
\end{equation}
The quotient is cocentral, meaning
\begin{equation}\label{eq:cocentral}
(\pi\otimes\id)\Delta(h)=(\pi\otimes\id)\Delta^{\cop}(h)
\qquad(h\in H).
\end{equation}
Furthermore, $K=H^{\co\pi}$ is a normal Hopf subalgebra and
$R\in K\otimes K$.
\end{lemma}

\begin{proof}
Let $T=(\pi\otimes\id)(R)\in B\otimes H$. The map
\[
f_R:H^{*\cop}\longrightarrow H,\qquad
f_R(p)=(p\otimes\id)(R)
\]
is a Hopf map. Thus
$L=f_R(\pi^*(B^*))$ is a Hopf subalgebra and $\dim (L)\mid\dim (B)$.
Its elements are the right coefficients of $T$.

Using $(\id\otimes\Delta)(R)=R_{13}R_{12}$, and applying
$\pi\otimes\id\otimes\pi$, gives
\[
(\id_B\otimes\id_H\otimes\pi)
(\id_B\otimes\Delta)(T)=T_{12}.
\]
Indeed, the image of $R_{13}$ is $1$ by the first hypothesis.
Taking coefficients in the first tensor factor proves
$L\subseteq H^{\co\pi}$.
By \eqref{eq:coinvariant-div}, $\dim (L)$ also divides
$\dim (H)/\dim( B)$, therefore   $L=\C$.
Applying the counit to the second factor of $T$ now gives $T=1\otimes1$.

The second hypothesis then yields $
1\otimes 1=(\pi\otimes\id)(R_{21}R)
  =(\pi\otimes\id)(R_{21})$,
whose tensor flip proves the other equality in
\eqref{eq:mixed-trivial}. Applying $\pi\otimes\id$ to the last
identity in \eqref{eq:qt} proves \eqref{eq:cocentral}. To explain normality, dualize \eqref{eq:cocentral}. For
$b\in B^*$ and $\varphi\in H^*$, evaluation on $h\in H$ gives
\[
(\pi^*(b)\varphi)(h)
 =b(\pi(h_{(1)}))\varphi(h_{(2)})
 =\varphi(h_{(1)})b(\pi(h_{(2)}))
 =(\varphi\pi^*(b))(h).
\]
Thus $\pi^*(B^*)$ is a central, hence normal, Hopf subalgebra of $H^*$.
Dualizing its exact sequence shows that $K=H^{\co\pi}$ is a normal
Hopf subalgebra and $
\C\longrightarrow K\longrightarrow H
\xrightarrow{\pi}B\longrightarrow\C
$ is exact.

Finally, \eqref{eq:mixed-trivial} and the two tensor identities in
\eqref{eq:qt} give
\[
(\id\otimes\id\otimes\pi)(\id\otimes\Delta)(R)=R_{12},
\qquad
(\id\otimes\pi\otimes\id)(\Delta\otimes\id)(R)=R_{13}.
\]
The first identity puts the right coefficients of $R$ in $K$, and the
second puts its left coefficients in $K$. Hence $R\in K\otimes K$.
\end{proof}

The following lemma is easy, and will play a role in the cocommutativity argument.
\begin{lemma}\label{lem:pointed}
Let $\cC$ be an integral braided fusion category over $\C$ with
square-free Frobenius--Perron dimension. If $\cC'$ is $\Vect$ or
$\sVec$, then $\cC$ is pointed.
\end{lemma}

\begin{proof}
An integral fusion category has a canonical
spherical structure with
$\dim(X)=\FPdim(X)$ for every simple object $X$
\cite[Corollary~9.6.6]{EGNO}, where $\dim(X)$ is the dimension determined by the spherical structure of $\cC$.
If $\cC'=\Vect$, then $\cD$ is modular; if $\cC'=\sVec$, then $\cC$ is super-modular. Then it follows from \cite[Proposition~8.14.6]{EGNO} and \cite[Corollary~3.4]{Yu} that
$\frac{\dim(\cC)}{\dim(X)^2}$
is an algebraic integer, which implies $\dim(X)=1$ as $\dim(\cC)$ is square-free.
\end{proof}
A symmetric fusion category $\cC$ is Tannakian if $\cC\cong\Rep(G)$ as braided fusion category \cite{DGNO,EGNO}, where $G$ is a finite group and the braiding of $\Rep(G)$ is given by the reflection of vector spaces.

\begin{theorem}\label{thm:cocommutative}
Every semisimple quasi-triangular Hopf algebra $(H,R)$ over $\C$ of square-free
dimension is cocommutative.
\end{theorem}
\begin{proof}
We adapt the method of  \cite[Lemma 5.2]{Natale} from a categorical view. Let $d=\dim (H)$ and $\cC=\Rep(H)$ with braiding $c=\tau\circ R$.
Let $\cE=\Rep(G)\subseteq\cC'$ be a  maximal Tannakian subcategory and $\cD:=\cC_G$. Then it follows from \cite[Proposition 4.30]{DGNO} that $
\cD'\simeq\Vect\ \text{or}\ \sVec$  and $\FPdim(\cD)=N/|G|$,  so $\cD$ is  pointed
by Lemma~\ref{lem:pointed}.

The fusion subcategory $\cE\subseteq\Rep (H)$ corresponds to a Hopf
quotient $\pi:H\twoheadrightarrow B$, of dimension $|G|$.
The induced triangular structure on $B$ is Tannakian, then
$(B,R_B)\cong(\C [G],1\otimes1)$ by \cite[Lemma 5.1]{Natale}.
Lemma~\ref{lem:cocentral} then shows that there is an exact sequence of Hopf algebras
\begin{equation}\label{eq:hopf-extension}
\C\longrightarrow K\longrightarrow H
\xrightarrow{\pi}\C [G]\longrightarrow\C,
\quad K=H^{\co\pi},\quad R\in K\otimes K,
\end{equation}
with $\pi$ cocentral. Therefore,  $\Rep(H)$ is a $G$-equivariantization
of $\Rep(K)$, with restriction to $K$ corresponding to the
forgetful functor \cite[Section~4.3]{Natale11}.
Consequently,
$\Rep(K)\cong\cC_G=\cD$
as fusion  categories.

Since $K$ is semisimple and $\Rep(K)$ is pointed. Thus $K$ is commutative, so
$K\cong\C^A$ for a finite group $A$. The restriction of $R$ makes
$\Rep(K)$ braided, hence  $A$  must be cyclic as its order is square-free.
The sequence \eqref{eq:hopf-extension} is therefore a cocentral extension
\[
\C\longrightarrow\C^A\longrightarrow H
\longrightarrow\C[G]\longrightarrow\C
\]
with $A$ cyclic.
Thus $H$ is cocommutative by \cite[Lemma 5.8]{Natale11}.
\end{proof}

Recall that a fiber functor  $F:\cC\to\Vect$ is an exact faithful tensor functor \cite{EGNO}.
\begin{corollary}\label{cor:category}
Let $\cC$ be a braided finite tensor category over $\C$. Suppose that
$\FPdim(\cC)=d$ is square-free and that there exists a fiber functor
$F\colon\cC\longrightarrow\Vect$.
Then $\cC$ is a fusion category and there is a tensor equivalence $
\cC\simeq\Rep(G)$
for a finite group $G$ of order $d$.
\end{corollary}

\begin{proof}The Tannakian reconstruction \cite[Theorem 5.3.12]{EGNO} show that there is a  finite-dimensional Hopf algebra $H$ and a
tensor equivalence $
\cC\cong\Rep(H)$
identifying $F$ with the forgetful functor. Transporting the braiding
along this equivalence equips $H$ with an $R$-matrix.
Thus $H$ is quasi-triangular and has square-free dimension.
Theorem~\ref{thm:cocommutative} shows that $H\cong\C [G]$ as a Hopf algebra, and then $\cC\cong\Rep (G)$ as a fusion category.
\end{proof}

\begin{remark}
The conclusion in Corollary~\ref{cor:category} concerns the underlying Hopf algebra. It does not assert
that $R=1\otimes1$, or even that the given braiding on $\Rep(H)$ is
symmetric. Therefore, the equivalence is
just a tensor equivalence; the braiding transported to $\Rep(G)$ may differ
from its usual symmetric braiding.
\end{remark}

\section*{Acknowledgments} Z. Yu  is supported by NSFC (no.\ 12571041) and Qinglan Project of Yangzhou University.

\author{			Jiacheng Bao, \thanks{Email:\,MX120240308@stu.yzu.edu.cn},\\
{Zhiqiang Yu, \thanks{Email:\,zhiqyumath@yzu.edu.cn}\\{Ying Zheng, \thanks{Email:\,yzhengmath@yzu.edu.cn}, \\
{\small School  of Mathematics,  Yangzhou University, Yangzhou 225002, China}}
\end{document}